\documentclass[12pt]{amsart}
\numberwithin{equation}{section}
\newtheorem{theorem}{Theorem}
\numberwithin{theorem}{section}

\numberwithin{lemma}{section}

\numberwithin{definition}{section}
\newtheorem{corollary}{Corollary}
\numberwithin{corollary}{section}

\numberwithin{proposition}{section}
\def\b{\begin{equation}}
\def\e{\end{equation}}
\newcommand{\ignore}[1]{}
\normalsize
\date {February 10, 2005}
\begin{document}
\pagenumbering{arabic} \pagenumbering{arabic}\setcounter{page}{1}
\tracingpages 1
\title{Sharp Hardy type inequalities on the Carnot Group}
\author{Ismail Kombe }
\address{Ismail Kombe, Mathematics Department\\ Dawson-Loeffler Science
\&Mathematics Bldg\\
Oklahoma City University \\
2501 N. Blackwelder, Oklahoma City, OK 73106-1493}
\email{ikombe@okcu.edu}
\begin{abstract}
In this paper we establish sharp weighted Hardy type inequalities
on the Carnot group with homogeneous dimension $Q\ge 3$.
\end{abstract}
\maketitle
\section{Introduction}

The classical Hardy inequality states that for $n\ge 3$
\begin{equation}
\int_{\mathbb{R}^n}|\nabla\phi(x)|^2dx\ge
\Big(\frac{n-2}{2}\Big)^2\int_{\mathbb{R}^n}
\frac{|\phi(x)|^2}{|x|^2}dx
\end{equation}
where $\phi\in C_0^{\infty}( \mathbb{R}^n\setminus \{0\})$ and the
constant $(\frac{n-2}{2})^2$ is sharp. There exists a large
literature dealing with the Hardy type inequalities on the
Euclidean space and, in particular, sharp inequalities which have
attracted considerable attention because of their application to
certain singular problems .

In this paper we investigate the existence  and the explicit
determination of constants $C$ and weight $q(x)$ on the Carnot
group $\mathbb{G}$ such that the Hardy type inequality

\begin{equation}
\int_{ \mathbb{G}}w(x)|\nabla_{ \mathbb{G}} \phi(x)|^2dx \ge
C\int_{ \mathbb{G}} q(x)|\phi(x)|^2dx
\end{equation}
holds for all $\phi\in C_0^{\infty}( \mathbb{R}^n\setminus
\{0\})$. Here  we consider a special weight function $w(x)$ which
is related to the fundamental solution of sub-Laplacian $ \Delta_{
\mathbb{G}}$  on the Carnot group $\mathbb{G}$.

It is well known that the Euclidean space $ \mathbb{R}^n$ with its
usual Abelian group structure is a trivial Carnot group and the
weighted Hardy type inequalities have been studied extensively. We
are concerned with the inequality (1.2) on the non-trivial Carnot
groups.

The simplest non trivial example of the Carnot group is given by
the Heisenberg group $ \mathbb{H}^n$. The following Hardy type
inequality on the Heisenberg group $ \mathbb{H}^n$ was first
proved by Garofalo and Lanconelli \cite {10}
\begin{equation}
\int_{ \mathbb{H}^n}|\nabla_{ \mathbb{H}^n}\phi|^2dzdt \ge
\Big(\frac{Q-2}{2}\Big)^2\int_{ \mathbb{H}^n}
(\frac{|z|^2}{|z|^4+t^2})\phi^2 dzdt
\end{equation}
where $\phi\in C_0^{\infty}(\mathbb{H}^n\setminus \{0\})$,
$Q=2n+2$ and the constant $(\frac{Q-2}{2})^2$ is sharp.

The Hardy type inequalities  on the Heisenberg group $
\mathbb{H}^n$ have also received considerable attention in recent
years. The $L^p$ version of the inequality (1.3) has been obtained
by Niu, Zhang, and Wang \cite {14}. A different proof of (1.3)
with the sharp constant $(\frac{Q-2}{2})^2$ has been given by
Goldstein and Zhang \cite {11}. In \cite {2}, D'Ambrosio obtained
weighted Hardy type inequalities on the Heisenberg group $
\mathbb{H}^n$.

Although we prove Hardy type inequality on the Carnot group with
an arbitrary step, we first establish a sharp weighted Hardy type
inequality on the Heisenberg group and extend this result to the
\textit{H}-type groups. The main point is that the fundamental
solution of the sub-Laplacian on the Heisenberg group and
\textit{H}-type groups are known explicitly (see section 3) but
not for the general  Carnot group. The proof of our theorem on the
Carnot group differs slightly in some steps from the the
Heisenberg group and \textit{H}-type group cases. The method we
apply here, inspired by the work of Allegretto \cite {1}, can be
applied to the Baouendi-Grushin type vector fields in that they do
not arise from any Carnot group.
\medskip

In order to state our theorems on the Carnot group $\mathbb{G}$,
we first recall the basic properties of Carnot group $\mathbb{G}$
and some well known results that will be used in the sequel. The
following section is largely taken from \cite {3}, \cite{7},
\cite{8}, \cite{9}, \cite {15} and \cite {16}.
\section {Carnot group}

 A Carnot group is a connected,
simply connected, nilpotent Lie group
$\mathbb{G}\equiv(\mathbb{R}^n,\cdot)$ whose Lie algebra
$\mathcal{G}$ admits a stratification. That is, there exist linear
subspaces $V_1, ..., V_k$  of $\mathcal{G}$ such that
\begin{equation} \mathcal{G}=V_1\oplus...\oplus V_k, \quad [V_1, V_i]=V_{i+1},
\quad \text{for}\quad i=1,2, ..., k-1 \quad\text {and}\quad [V_1,
V_k]=0
\end{equation}
 where $[V_1, V_i]$ is the subspace  of $\mathcal{G}$
generated by the elements $[X,Y]$ with $X\in V_1$ and $Y\in V_i$.
This defines a $k$-step Carnot group and integer $k\ge 1$ which is
called the step of $ \mathbb{G}$.

Via the exponential map, it is possible to induce on $\mathbb{G}$
a family of automorphisms of the group, called dilations,
$\delta_{\lambda}: \mathbb{R}^n\longrightarrow \mathbb{R}^n
(\lambda>0)$ such that
\[\delta_{\lambda}(x_1, ...,x_n)=(\lambda^{\alpha_1} x_1,...,\lambda^{\alpha_n}x_n)\]
where $1=\alpha_1=...=\alpha_m<\alpha_{m+1}\le ...\le \alpha_n$
are integers and $m=\text{dim}(V_1)$.

The group law can be written in the following form
\begin{equation}
x\cdot y=x+y+P(x,y), \quad x, y\in \mathbb{R}^n
\end{equation}
where $ P:\mathbb{R}^n\times \mathbb{R}^n\longrightarrow
\mathbb{R}^n$ has polynomial components and $P_1=...=P_m=0$. Note
that the inverse $x^{-1}$ of an element $x\in \mathbb{G}$ has the
form $x^{-1}=-x=(-x_1,....-x_n)$.

Let $X_1, . . ., X_m$ be a family of left invariant vector fields
that is an orthonormal basis of $ V_1\equiv\mathbb{R}^m$ at the
origin, that is, $X_1(0)=\partial_{x_1}, . . ., X_m(0)=\partial
_{x_m}$. The vector fields $X_j$  have polynomial coefficients and
can be assumed to be of the form
\[X_j(x)=\partial_j+\sum_{i=j+1}^n a_{ij}(x)\partial_i, \quad
X_j(0)=\partial_j, j=1, . . ., m,\] where  each polynomial
$a_{ij}$ is homogeneous with respect to the dilations of the
group, that is
$a_{ij}(\delta_{\lambda}(x))=\lambda^{\alpha_i-\alpha_j}a_{ij}(x)$.
The horizontal  gradient on the Carnot group $\mathbb{G}$ is the
vector valued operator
\[\nabla_{\mathbb{G}}=(X_1, . . ., X_m)\] where $X_1, . . ., X_m$ are the generators of $ \mathbb{G}$.
The sub-Laplacian is the second-order partial differential
operator on $\mathbb{G}$  given by
\[\Delta_{ \mathbb{G}}=\sum_{j=1}^m X_j^2.\]
The fundamental solution $u$ for $\Delta_{\mathbb{G}}$ is defined
to be a weak solution to the equation
\[-\Delta_{\mathbb{G}}u=\delta\] where $\delta$ denotes the Dirac distribution with singularity at
the neutral element $0$ of $ \mathbb{G}$. In \cite{7} Folland
proved that in any Carnot group $ \mathbb{G}$,  there exists a
homogeneous norm $N$ such that
\[u=N^{2-Q}\] is a fundamental solution for $\Delta_{
\mathbb{G}}$ ( see also \cite{4}).

We now set $N(x):=u^{\frac{1}{2-Q}}$ if $x\neq 0$ and $N(0):=0$.
We recall that a homogeneous norm on $ \mathbb{G}$ is a continuous
 function $N:\mathbb{G}\longrightarrow [0, \infty )$ smooth away from the origin
which satisfies the conditions : $N(\delta_{\lambda}(x))=\lambda
N(x)$, $N(x^{-1})=N(x)$ and $N(x)=0$ iff $x=0$.

 The curve $\gamma:[a,b]\subset \mathbb{R}\longrightarrow
\mathbb{G}$ is called horizontal if its tangents lie in $V_1$,
i.e,  $\gamma'(t)\in \text{\textit{span}}\{X_1, . . ., X_m\}$ for
all $t$. Then, the Carnot-Car\'ethedory distance $d_{CC}(x,y)$
between two points $x,y\in \mathbb{G}$ is defined to be the
infimum of all horizontal lengths $\int_a^b \langle \gamma'(t),
\gamma'(t)\rangle^{1/2} dt$ over all horizontal curves
$\gamma:[a,b]\longrightarrow \mathbb{G}$ such that $\gamma(a)=x$
and $\gamma(y)=b$. Notice that $d_{CC}$ is a homogeneous norm and
satisfies the invariance property

\[d_{CC}(z\cdot x, z\cdot y)=d_{CC}(x,y), \quad \forall\, x,y,
z\in \mathbb{G},\] and is homogeneous of degree one with respect
to the dilation $\delta_{\lambda}$, i.e.
\[d_{CC}(\delta_{\lambda}(x), \delta_{\lambda}(y))=\lambda
d_{CC}(x,y), \quad \forall\, x, y, z \in \mathbb{G}, \forall\,
\lambda>0.\]

The Carnot-Careth\'edory balls are defined by $B(x, R)=\{y\in
\mathbb{G} |d_{CC}(x,y)<R\}$. By left-translation and dilation, it
is easy to see that the Haar measure of $B(x, R)$ is proportional
by $R^Q$. More precisely
\[ |B(x, R)|=R^Q |B(x,1)|=R^Q |B(0,1)|\] where \[Q=\sum_{j=1}^k j (\text{dim}V_j)\] is the
homogeneous dimension of $ \mathbb{G}$.

\section{Hardy type inequality on the Carnot group of step 2}
Among the Carnot groups of step two,  the Heisenberg group  and
Heisenberg type (\textit{H}-type) groups are of particular
significance. These groups appear naturally in analysis, geometry,
representation theory and mathematical physics. In this section,
we first prove Hardy type inequality on the Heisenberg group and
we extend this result to the \textit{H}-type group.

\noindent {\textbf{Heisenberg group}.} The Heisenberg group $
\mathbb{H}^n$ is an example of a noncommutative Carnot group.
Denoting points in $ \mathbb{H}^n$ by $(z,t)$ with $z=(z_1, . . .,
z_n)\in \mathbb{C}^n$ and $t\in \mathbb{R}$ we have the group law
given as
\[(z,t)\circ (z',t')=(z+z', t+t'+2\sum_{j=1}^n
Im(z_j\bar{z}_j'))\] With the notation $z_j=x_j+iy_j$, the
horizontal space $V_1$ is spanned by the basis
\[X_j=\frac{\partial}{\partial x_j}
+2y_j\frac{\partial }{\partial t}\quad\text{and} \quad
Y_j=\frac{\partial}{\partial y_j}-2x_j\frac{\partial }{\partial
t}.\] The one dimensional center $V_2$ is spanned by the vector
field $T=\frac{\partial}{\partial t}$.  We have the commutator
relations $[X_j, Y_j]=-4T$, and all other brackets of $\{X_1, Y_1,
. . ., X_n, Y_n\}$ are zero. The sub-elliptic gradient is the $2n$
dimensional vector field given by
\[\nabla_{ \mathbb{H}^n}=(X_1, . . ., X_n, Y_1, . . ., Y_n)\] and
the Kohn Laplacian on $ \mathbb{H}^n$ is the operator \[\Delta_{
\mathbb{H}^n}=\sum_{j=1}^n (X_j^2+Y_j^2).\] A homogeneous norm on
$ \mathbb{H}^n$ is given by
\[\rho=|(z,t)|=(|z|^4+t^2)^{1/4}\] and the homogeneous dimension of $\mathbb{H}^n$ is $Q=2n+2$.

A remarkable analogy between Kohn Laplacian and the classical
Laplace operator has been obtained by Folland \cite{6}. He found
that the fundamental solution of $-\Delta_{ \mathbb{H}^n}$ with
pole zero is given by
\[\Psi(z,t)=\frac{c_Q}{\rho(z,t)^{Q-2}}\quad\text{where}\quad
c_Q=\frac{2^{(Q-2)/2}\Gamma((Q-2)/4)^2}{\pi^{Q/2}}.\]

We now prove the following theorem  on the Heisenberg group $
\mathbb{H}^n$ (See \cite {2} for a different proof). In the
various integral inequalities below (Section 3 and Section 4), we
allow the values of the integrals on the left-hand sides to be
$+\infty$.

\begin{theorem}
Let $\alpha\in \mathbb{R}$ and $\phi\in
C_0^{\infty}(\mathbb{H}^n\setminus \{0,0\})$. Then we have :
\[\int_{
\mathbb{H}^n}\rho^{\alpha}|\nabla_{\mathbb{H}^n}\phi|^2dzdt\ge
\Big(\frac{Q+\alpha-2}{2}\Big)^2\int_
{\mathbb{H}^n}\rho^{\alpha}\frac{|z|^2}{\rho^4}\phi^2dzdt\] where
$\rho=(|z|^4+l^2)^{1/4}$ is the homogeneous norm on $
\mathbb{H}^n$. Moreover, the constant $(\frac{Q+\alpha-2}{2})^2$
is sharp.
\end{theorem}
\proof Let $\phi=\rho^{\beta}\psi$ where $\beta\in
\mathbb{R}\setminus \{0\}$ and $\psi \in
C_0^{\infty}(\mathbb{H}^n\setminus \{0,0\})$. A direct calculation
shows that
\begin{equation}
\rho^{\alpha}|\nabla_{ \mathbb{H}^n}\phi|^2
=\beta^2\rho^{\alpha+2\beta-2}|\nabla_{\mathbb{H}^n}\rho|^2\psi^2+
2\beta\rho^{\alpha+2\beta-1}\psi\nabla_{
\mathbb{H}^n}\rho\cdot\nabla_{
\mathbb{H}^n}\psi+\rho^{\alpha+2\beta}|\nabla_{\mathbb{H}}\psi|^2.
\end{equation}
It is easy to see that
\[|\nabla_{ \mathbb{H}^n}\rho|^2=\frac{|z|^2}{\rho^2}\]
and integrating (3.1) over $ \mathbb{H}^n$, we get
\begin{equation}
\begin{aligned}\int_{ \mathbb{H}^n}\rho^{\alpha}|\nabla_{
\mathbb{H}^n}\phi|^2dzdt&=\int_{
\mathbb{H}^n}\beta^2\rho^{\alpha+2\beta-4}|z|^2\psi^2dzdt+\int_{
\mathbb{H}^n} 2\beta\rho^{\alpha+2\beta-1}\psi\nabla_{
\mathbb{H}^n}\rho\cdot\nabla_{ \mathbb{H}^n}\psi dzdt\\&+\int_{
\mathbb{H}^n}\rho^{\alpha+2\beta}|\nabla_{\mathbb{H}^n}\psi|^2dzdt\\
\end{aligned}
\end{equation}
Applying integration by parts to the middle integral on the
right-hand side of (3.2), we obtain

\begin{equation}
\begin{aligned}
\int_{ \mathbb{H}^n}\rho^{\alpha}|\nabla_{
\mathbb{H}^n}\phi|^2dzdt&=\int_{
\mathbb{H}^n}\beta^2\rho^{\alpha+2\beta-4}|z|^2\psi^2dzdt-\frac{\beta}{\alpha+2\beta}
\int_{ \mathbb{H}^n}\Delta_{
\mathbb{H}^n}(\rho^{\alpha+2\beta})\psi^2dzdt\\&+\int_{
\mathbb{H}^n}\rho^{\alpha+2\beta}|\nabla_{\mathbb{H}^n}\psi|^2dzdt.
\end{aligned}
\end{equation}
One can show that
\begin{equation}\Delta_{
\mathbb{H}^n}(\rho^{\alpha+2\beta})=|z|^2\rho^{\alpha+2\beta-4}(\alpha+2\beta)(\alpha+2\beta+Q-2).
\end{equation}
Substituting (3.4) into (3.3) gives the following

\[\begin{aligned} \int_{
\mathbb{H}^n}\rho^{\alpha}|\nabla_{ \mathbb{H}^n}\phi|^2dzdt &=
(\beta^2-\beta(\alpha+2\beta+Q-2))\int_{
\mathbb{H}^n}\rho^{\alpha+2\beta-4}|z|^2\psi^2dzdt+ \int_{
\mathbb{H}^n}
\rho^{\alpha+2\beta}|\nabla_{\mathbb{H}^n}\psi|^2dzdt\\
&\ge(-\beta^2-\beta(\alpha+Q-2))\int_{
\mathbb{H}^n}\rho^{\alpha+2\beta-4}|z|^2\psi^2dzdt.
\end{aligned}
\]
We now choose  $\beta=\frac{2-\alpha-Q}{2}$ ( Note that the
quadratic equation  $-\beta^2-\beta(\alpha+Q-2)$ reaches its
maximum value at  $\beta=\frac{2-\alpha-Q}{2}$) and noting that
$\psi=\rho^{-\beta}\phi$, we have the following inequality
\begin{equation}
\int_{ \mathbb{H}^n}\rho^{\alpha}|\nabla_{
\mathbb{H}^n}\phi|^2dzdt \ge
\Big(\frac{Q+\alpha-2}{2}\Big)^2\int_{
\mathbb{H}^n}\rho^{\alpha}\frac{|z|^2}{\rho^4}\phi^2dzdt.
\end{equation}
\endproof
\medskip

\noindent{\textbf{Heisenberg type group}.} Another important model
of Carnot groups are the \textit{H}-type (Heisenberg type) groups
which were introduced by Kaplan \cite {13} as direct
generalizations of the Heisenberg group $\mathbb{H}^n$. An
\textit{H}-type group is a Carnot group with a two-step Lie
algebra $ \mathcal{G}=V_1\oplus V_2$ and an inner product
$\langle, \rangle$ in $ \mathcal{G}$ such that the linear map
\[ J:V_2\longrightarrow \text{End}V_1,\] defined by the condition
\[\langle J_z(u), v\rangle=\langle z, [u,v]\rangle,\quad u,v\in V_1, z\in V_2\] satisfies
\[J_z^2=-||z||^2\mathbf{Id}\]
for all $z\in V_2$, where $||z||^2= \langle z,z\rangle$.

Sub-Laplacian is defined in terms of a fixed basis $X_1,. . . ,
X_m$ for $V_1$:

\begin{equation}
\Delta_\mathbb{G}=\sum_{i=1}^mX_i^2.
\end{equation}

The exponential mapping of a simply connected Lie group is an
analytic diffeomorphism. One can then define analytic mappings
$v:\mathbb{G} \longrightarrow V_1$ and
$z:\mathbb{G}\longrightarrow V_2$ by
\[x=\text{exp}(v(x)+z(x))\] for every $x\in \mathbb{G}$.
In \cite {13} Kaplan proved that there exists a constant $c>0$
such that the function
\[\Phi (x)=c\Big(|v(x)|^4+16|z(x)|^2\Big)^{\frac{2-Q}{4}}\]
is a  fundamental solution for the operator $ -\Delta_\mathbb{G}$.
We note that
\begin{equation}
K(x)=\Big(|v(x)|^4+16|z(x)|^2\Big)^{\frac{1}{4}}
\end{equation}
defines a homogeneous norm  and $Q=m+2k$ is the homogeneous
dimension of $\mathbb{G}$ where $m=$dim$V_1$ and $k=$dim$V_2$.
This result generalized Folland's fundamental solution for the
Heisenberg group $ \mathbb{H}^n$ \cite{6}. Note that Kaplan's
results provides us an explicit fundamental solution for $
\Delta_\mathbb{G}$ on \textit{H}-type groups. We should also
mention that the explicit fundamental solutions for the
sub-elliptic $p$-Laplacian on the \textit{H}-type groups were
obtained by Capogna, Danielli and Garofalo \cite{5},  Heinonen and
Holopainen \cite{12}.

We cite, without proof of the following, useful formulas which can
be found in \cite{5}: Let $u$ be a radial function, i.e.,
$u(x)=f(K(x))$ where $f\in C(\mathbb{R})$ then

\[ |\nabla_{ \mathbb{G}}u|^2=\frac{|v|^2}{K^2}|f'(K)|^2.\] Moreover  if  $u(x)=f(K(x))$  and $f\in
C^2(\mathbb{R})$ then
\begin{equation}
\begin{aligned}
\Delta_{\mathbb{G}}u&=|\nabla_{\mathbb{G}}K(x)|^2\Big[f''(K)+\frac{Q-1}{K}f'(K)\Big]\\
&=\frac{v^2}{K^2} \Big[f''(K)+\frac{Q-1}{K}f'(K)\Big]
\end{aligned}
\end{equation}
 at every point $x\in \mathbb{G}\setminus \{0\}$ where
$f'(K(x))\neq 0$.
\medskip

We now have the following theorem on the \textit{H}-type group :

\begin{theorem}
Let $ \mathbb{G}$ be an \textit{H}-type group with homogeneous
dimension $Q=m+2k$ and let $\alpha\in \mathbb{R}$ and $\phi\in
C_0^{\infty}(\mathbb{G}\setminus \{0\})$. Then the following
inequality is valid :
\begin{equation}\int_{
\mathbb{G}}K^{\alpha}|\nabla_{\mathbb{G}}\phi|^2dx\ge
\Big(\frac{Q+\alpha-2}{2}\Big)^2\int_
{\mathbb{G}}K^{\alpha}\frac{|v|^2}{K^4}\phi^2dx\end{equation}
where $K(x)=(|v(x)|^4+16|z(x)|^2)^{1/4}$. Moreover, the constant
$(\frac{Q+\alpha-2}{2})^2$ is sharp.
\end{theorem}
\proof The proof is identical to the Heisenberg group case.
\endproof

\section{Hardy-type inequality on the Carnot group of arbitrary step }
In this section,  we consider the Carnot group $ \mathbb{G}$ of
any step $k$ with the homogeneous norm $N=u^{1/(2-Q)}$ associated
to Folland's solution $u$ for the sub-Laplacian $\Delta_{
\mathbb{G}}$ \cite{7}. We have the following theorem:
\begin{theorem} Let $\mathbb{G}$ be a Carnot group with homogeneous dimension $Q\ge 3$  and  let $\phi\in
C_0^{\infty}(\mathbb{G}\setminus \{0\})$, $\alpha\in \mathbb{R}$,
$Q+\alpha-2>0$. Then the following inequality is valid
\begin{equation}\int_{\mathbb{G}} N^{\alpha}|\nabla_{\mathbb{G}} \phi|^2dx
\ge \Big(\frac{Q+\alpha-2}{2}\Big)^2 \int_{\mathbb{G}}
N^{\alpha}\frac{|\nabla_{\mathbb{G}} N|^2}{N^2}\phi
^2dx.\end{equation} Here $N=u^{1/(2-Q)}$ is the homogeneous norm
associated with the fundamental solution $u$ for the sub-Laplacian
$\Delta_{ \mathbb{G}}$. Furthermore,  the constant
$C(Q,\alpha)=(\frac{Q+\alpha-2}{2})^2$ is sharp.
\end{theorem}
\proof

 Let $\phi=N^{\beta}\psi$  where $\psi\in
C_0^{\infty}(\mathbb{G}\setminus \{0\})$ and $\beta\in
\mathbb{R}\setminus \{0\}$. A direct calculation
 shows that
 \begin{equation}
|\nabla_{\mathbb{G}}(N^{\beta}\psi)|^2=\beta^2N^{2\beta-2}|\nabla_{\mathbb{G}}
N|^2\psi^2+2\beta N^{2\beta-1}\psi\nabla_{\mathbb{G}} N\cdot
\nabla_{\mathbb{G}} \psi+N^{2\beta}|\nabla_{\mathbb{G}}\psi|^2.
\end{equation}
Multiplying both sides of (4.2) by the $N^{\alpha}$ and applying
integration by parts over $ \mathbb{G}$ gives

\begin{equation}\begin{aligned}
\int_{\mathbb{G}}N^{\alpha}|\nabla_{\mathbb{G}}\phi|^2dx
&=\beta^2\int_{\mathbb{G}} N^{\alpha+2\beta-2}|\nabla_{\mathbb{G}}
N|^2\psi ^2dx-\frac{\beta}{\alpha+2\beta}\int_{\mathbb{G}}
\Delta_{\mathbb{G}}(N^{\alpha+2\beta})\psi^2dx\\&+\int_{\mathbb{G}}N^{\alpha+2\beta}|\nabla_{\mathbb{G}}\psi|^2dx\\
& \ge \beta^2\int_{\mathbb{G}}
N^{\alpha+2\beta-2}|\nabla_{\mathbb{G}} N|^2\psi
^2dx-\frac{\beta}{\alpha+2\beta}\int_{\mathbb{G}}
\Delta_{\mathbb{G}}(N^{\alpha+2\beta})\psi^2dx.
\end{aligned}
\end{equation}
A straightforward calculation shows that

\begin{equation}
-\frac{\beta}{\alpha+2\beta} \Delta_{
\mathbb{G}}(N^{\alpha+2\beta})=-\beta(\alpha+2\beta+Q-2)N^{\alpha+2\beta-2}|\nabla_{\mathbb{G}}N|^2
-\frac{\beta}{2-Q}N^{\alpha+2\beta+Q-2}\Delta_{\mathbb{G}}u.
\end{equation}
Substituting (4.4) into (4.3) and  using the fact that
$\psi^2=N^{-2\beta}\phi^2$, we get the following :

\[
\int_{\mathbb{G}}N^{\alpha}|\nabla_{\mathbb{G}}\phi|^2dx\ge (
-\beta^2-\beta(\alpha+Q-2)\int_{
\mathbb{G}}N^{\alpha}\frac{|\nabla_{ \mathbb{G}}
N|^2}{N^2}\phi^2dx-\frac{\beta}{2-Q}\int_{
\mathbb{G}}(\Delta_{\mathbb{G}}u)N^{\alpha+Q-2}\phi^2dx.
\]
Since $u$ is the fundamental solution of sub-Laplacian
$\Delta_{\mathbb{G}}$ on the Carnot group $\mathbb{G}$, we get
\[\int_{
\mathbb{G}}(\Delta_{\mathbb{G}}u)N^{\alpha+Q-2}\phi^2dx=
N^{\alpha+Q-2}(0)\phi^2(0)=0.\]  We now obtain

\[\int_{\mathbb{G}}N^{\alpha}|\nabla_{\mathbb{G}}\phi|^2dx\ge (
-\beta^2-\beta(\alpha+Q-2)\int_{
\mathbb{G}}N^{\alpha}\frac{|\nabla_{ \mathbb{G}}
N|^2}{N^2}\phi^2dx.\] Choosing
\[\beta=\frac{2-Q-\alpha}{2}\]  gives the following sharp
inequality
\begin{equation}
\int_{\mathbb{G}} N^{\alpha}|\nabla_{\mathbb{G}} \phi|^2dx \ge
\Big(\frac{Q+\alpha-2}{2}\Big)^2 \int_{\mathbb{G}}
N^{\alpha}\frac{|\nabla_{\mathbb{G}} N|^2}{N^2}\phi ^2dx.
\end{equation}
\endproof

An immediate consequence of the Hardy type inequality (4.1) is the
following corollary, known as the \textit{uncertainty principle}.
\begin{corollary}(Uncertainty principle).
Let $\mathbb{G}$ be a Carnot group with homogeneous dimension
$Q\ge 3$. Then for every $\phi\in C_0^{\infty}(\mathbb{G}\setminus
\{0\})$
\begin{equation}
\Big(\int_{\mathbb{G}}
N^2|\nabla_{\mathbb{G}}N|^2\phi^2dx\Big)\Big(\int_{\mathbb{G}}|\nabla_{\mathbb{G}}
\phi|^2dx\Big)\ge
\Big(\frac{Q-2}{2}\Big)^2\Big(\int_{\mathbb{G}}|\nabla_{\mathbb{G}}N|^2\phi^2
dx\Big)^2.
\end{equation}
Here $N=u^{1/(2-Q)}$ is the homogeneous norm associated with the
fundamental solution $u$ for the sub-Laplacian $\Delta_{
\mathbb{G}}$.
\end{corollary}
\noindent\textbf{Remark}. In the Abelian case, when
$\mathbb{G}=\mathbb{R}^n$ with the ordinary dilations, one has
$\mathcal{G}=V_1=\mathbb{R}^n$ so that $Q=n$. Now it is clear that
the inequality (4.1) with the homogeneous norm $N(x)=|x|$ and
$\alpha=0$ recovers the Hardy's inequality (1.1).

\bibliographystyle{amsalpha}

\end{document}